\documentclass[leqno]{amsart}
\usepackage[normal]{boilerart}
\usepackage{mathrsfs}
\usepackage{tikz-cd}

\DeclareMathOperator{\Tot}{Tot}

\newcommand{\leftnat}[1]{\vphantom{#1}_{\natural}\mskip-1mu{#1}}

\def\reflectedop#1#2{\mathop{\reflectbox{$#1{#2}$}}}

\title{Fibrations in $\infty$-category theory}

\author{Clark Barwick}
\address{Department of Mathematics, Massachusetts Institute of Technology, 77 Massachusetts Avenue, Cambridge, MA 02139-4307, USA}
\email{clarkbar@math.mit.edu}

\author{Jay Shah}
\address{Department of Mathematics, Massachusetts Institute of Technology, 77 Massachusetts Avenue, Cambridge, MA 02139-4307, USA}
\email{jshah@math.mit.edu}

\thanks{We are thankful to Bob Bruner and Emily Riehl for their helpful comments on this paper.}

\begin{document}

\begin{abstract} In this short expository note, we discuss, with plenty of examples, the bestiary of fibrations in quasicategory theory. We underscore the simplicity and clarity of the constructions these fibrations make available to end-users of higher category theory.
\end{abstract}

\maketitle

\tableofcontents

The theory of $\infty$-categories -- as formalized in the model of quasicategories -- offers two ways of specifying homotopy theories and functors between them.

First, we may describe a homotopy theory via a homotopy-coherent universal property; this is a widely appreciated advantage, and it's a feature that any sufficiently well-developed model of $\infty$-categories would have. For example, the $\infty$-category $\Top$ of spaces is the free $\infty$-category generated under (homotopy) colimits by a single object \cite[Th. 5.1.5.6]{HTT}.

The second way of specifying $\infty$-categories seems to be less well-loved: this is the ability to perform completely explicit constructions with excellent formal properties. This allows one to avoid the intricate workarounds that many of us beleaguered homotopy theorists have been forced to deploy in order solve infinite hierarchies of homotopy coherence problems. This feature seems to be peculiar to the model of quasicategories, and the main instrument that makes these explicit constructions possible is the theory of \emph{fibrations} of various sorts. In this \'etude, we study eight sorts of fibrations of quasicategories in use today -- left, right, Kan, inner, iso (\textsc{aka} categorical), cocartesian, cartesian, and flat -- and we discuss the beautifully explicit constructions they provide.

In the end, $\infty$-category theory as practiced today combines these two assets, and the result is a powerful amalgam of universal characterizations and crashingly explicit constructions. We will here focus on the underappreciated latter feature, which provides incredibly concrete constructions to which we would not otherwise have access. We have thus written this under the assumption that readers are more or less familiar with the content of the first chapter of Lurie's book \cite{HTT}.

\section{Left, right, and Kan fibrations} The universal property of the $\infty$-category $\Top$ we offered above certainly characterizes it up to a contractible choice, but it doesn't provide any simple way to specify a functor \emph{into} $\Top$.

At first blush, this looks like very bad news: after all, even if $C$ is an ordinary category, to specify a functor of $\infty$-categories $F\colon\fromto{C}{\Top}$, one has to specify an extraordinary amount of information: one has to give, for every object $a\in C$, a space $F(a)$; for every morphism $f\colon\fromto{a}{b}$, a map $F(f)\colon\fromto{F(a)}{F(b)}$; for every pair of composable morphisms $f\colon\fromto{a}{b}$ and $g\colon\fromto{b}{c}$, a homotopy $F(gf)\simeq F(g)F(f)$; for every triple of composable morphisms, a homotopy of homotopies; etc., etc., \emph{ad infinitum}. 

However, ordinary category theory suggests a way out: suppose $F\colon\fromto{C}{\Set}$ a functor. One of the basic tricks of the trade in category theory is to build a category $\Tot F$, sometimes called the \emph{category of elements} of $F$. The objects of $\Tot F$ are pairs $(a,x)$, where $a\in C$ is an object and $x\in F(a)$ is an element; a morphism $\fromto{(a,x)}{(b,y)}$ of $\Tot F$ is a morphism $f\colon\fromto{a}{b}$ such that $F(f)(x)=y$.

The category $\Tot F$, along with the projection $p\colon\fromto{\Tot F}{C}$, is extremely useful for studying the functor $F$. For example, the set of sections of $p$ is a limit of $F$, and the set $\pi_0(\Tot F)$ of connected components is a colimit of $F$. In fact, the assignment $\goesto{F}{\Tot F}$ is an equivalence of categories between the category of functors $\fromto{C}{\Set}$ and those functors $\fromto{X}{C}$ such that for any morphism $f\colon\fromto{a}{b}$ of $C$ and for any object $x\in X$ with $p(x)=a$, there exists a unique morphism $\phi\colon\fromto{x}{y}$ with $p(\phi)=f$. (These functors are sometimes called \emph{discrete opfibrations}.) In other words, functors $\fromto{C}{\Set}$ correspond to functors $\fromto{X}{C}$ such that for any solid arrow commutative square
\begin{equation*}
\begin{tikzpicture}[baseline]
\matrix(m)[matrix of math nodes,
row sep=4ex, column sep=4ex,
text height=1.5ex, text depth=0.25ex]
{\Lambda^1_0 & NX \\
\Delta^1 & NC, \\ };
\path[>=stealth,->,font=\scriptsize]
(m-1-1) edge node[above]{} (m-1-2)
edge[right hook->] (m-2-1)
(m-1-2) edge node[right]{$p$} (m-2-2)
(m-2-1) edge node[below]{} (m-2-2)
(m-2-1) edge[dotted] (m-1-2);
\end{tikzpicture}
\end{equation*}
there exists a unique dotted lift.

We may therefore hope that, instead of working with functors from $C$ to $\Top$, one might work with suitable $\infty$-categories \emph{over} $C$ instead. To make this work, we need to formulate the $\infty$-categorical version of this condition. To this end, we adopt the same attitude that permits us to arrive at the definition of an $\infty$-category: instead of demanding a single unique horn filler, we demand a whole hierarchy of horn fillers, none of which we require to be unique. The hierarchy ensures that the filler at any stage is unique up to a homotopy that is unique up to a homotopy, etc., etc., \emph{ad infinitum}.
\begin{dfn} A \emph{left fibration} is a map $p\colon\fromto{X}{S}$ of simplicial sets such that for any integer $n\geq 1$ and any $0\leq k<n$ and any solid arrow commutative square
\begin{equation*}
\begin{tikzpicture}[baseline]
\matrix(m)[matrix of math nodes,
row sep=4ex, column sep=4ex,
text height=1.5ex, text depth=0.25ex]
{\Lambda^n_k & X \\
\Delta^n & S, \\ };
\path[>=stealth,->,font=\scriptsize]
(m-1-1) edge node[above]{} (m-1-2)
edge[right hook->] (m-2-1)
(m-1-2) edge node[right]{$p$} (m-2-2)
(m-2-1) edge node[below]{} (m-2-2)
(m-2-1) edge[dotted] (m-1-2);
\end{tikzpicture}
\end{equation*}
there exists a dotted lift.

Dually, a \emph{right fibration} is a map $p\colon\fromto{X}{S}$ of simplicial sets such that for any integer $n\geq 1$ and any $0<k\leq n$ and any solid arrow commutative square
\begin{equation*}
\begin{tikzpicture}[baseline]
\matrix(m)[matrix of math nodes,
row sep=4ex, column sep=4ex,
text height=1.5ex, text depth=0.25ex]
{\Lambda^n_k & X \\
\Delta^n & S, \\ };
\path[>=stealth,->,font=\scriptsize]
(m-1-1) edge node[above]{} (m-1-2)
edge[right hook->] (m-2-1)
(m-1-2) edge node[right]{$p$} (m-2-2)
(m-2-1) edge node[below]{} (m-2-2)
(m-2-1) edge[dotted] (m-1-2);
\end{tikzpicture}
\end{equation*}
there exists a dotted lift.

Of course \emph{Kan fibration} is a map of simplicial sets that is both a left and right fibration.
\end{dfn}

To understand these notions, we should begin with some special cases.
\begin{exm} For any simplicial set $X$, the unique map $\fromto{X}{\Delta^0}$ is a left fibration if and only if $X$ is an $\infty$-groupoid (i.e., a Kan complex). Indeed, we see immediately that $X$ is an $\infty$-category, so to conclude that $X$ is an $\infty$-groupoid, it suffices to observe that the homotopy category $hX$ is a groupoid; this follows readily from the lifting condition for the horn inclusion $\into{\Lambda^2_0}{\Delta^2}$. 
\end{exm}

Since pullbacks of left fibrations are again left fibrations, we conclude immediately that the fibers of a left fibration are $\infty$-groupoids.

\begin{exm}[\protect{\cite[Cor. 2.1.2.2]{HTT}}] If $C$ is an $\infty$-category and $x\in C_0$ is an object, then recall that one can form the \emph{undercategory} $C_{x/}$ uniquely via the following functorial bijection:
\[\Mor(K,C_{x/})\cong\Mor(\Delta^0\star K,C)\times_{\Mor(\Delta^0,C)}\{x\}.\]
The inclusion $\into{K}{\Delta^0\star K}$ induces a forgetful functor $p\colon\fromto{C_{x/}}{C}$. The key fact (due to Joyal) is that $p$ is a left fibration; in particular, $C_{x/}$ is an $\infty$-category \cite[Cor. 2.1.2.2]{HTT}.

The fiber of $p$ over a vertex $y\in C_0$ is the $\infty$-groupoid whose $n$-simplices are maps $f\colon\fromto{\Delta^{n+1}}{C}$ such that $f({\Delta^{\{0\}}})=x$ and $f|_{\Delta^{\{1,\dots,n+1\}}}$ is the constant map at $y$. In other words, it is $\Hom_C^L(x,y)$ in the notation of \cite[Rk. 1.2.2.5]{HTT}.

Dually, we can define the \emph{overcategory} $C_{/x}$ via
\[\Mor(K,C_{x/})\cong\Mor(K\star\Delta^0,C)\times_{\Mor(\Delta^0,C)}\{x\},\]
and the inclusion $\into{K}{K\star\Delta^0}$ induces a right fibration $q\colon\fromto{C_{/x}}{C}$. The fiber of $q$ over a vertex $y\in C_0$ is then $\Hom_C^R(y,x)$.
\end{exm}

\begin{subexm} In the introduction, we mentioned that $\Top$ is the $\infty$-category that is freely generated under colimits by a single object $\ast$. This generator turns out to be the terminal object in $\Top$. Let us write $\Top_{\ast}$ for the overcategory $\Top_{\ast/}$. The forgetful functor $\fromto{\Top_{\ast}}{\Top}$ is a left fibration.

The fiber over a vertex $X\in\Top_0$ is the $\infty$-groupoid $\Hom_{\Top}^L(\ast,X)$, which we will want to think of as a model for $X$ itself.
\end{subexm}

Here is the theorem that is going to make that possible:
\begin{thm}[Joyal] Suppose $C$ an $\infty$-category (or more generally, any simplicial set). For any functor $F\colon\fromto{C}{\Top}$, we may consider the left fibration
\[\fromto{\Top_{\ast}\times_{\Top,F}C}{C}.\]
This defines an equivalence of $\infty$-categories
\[\equivto{\Fun(C,\Top)}{\categ{LFib}(C)},\]
where $\categ{LFib}(C)$ is the simplicial nerve of the full simplicial subcategory of $s\Set_{/C}$ spanned by the left fibrations.

Dually, for any functor $G\colon\fromto{C}{\Top^{\op}}$, we may consider the right fibration
\[\fromto{\Top_{\ast}^{\op}\times_{\Top^{\op},G}C}{C}.\]
This defines an equivalence of $\infty$-categories
\[\equivto{\Fun(C^{\op},\Top)}{\categ{RFib}(C)},\]
where $\categ{RFib}(C)$ is the simplicial nerve of the full simplicial subcategory of $s\Set_{/C}$ spanned by the right fibrations.
\end{thm}

\noindent A left fibration $p\colon\fromto{X}{C}$ is said to be \emph{classified by} $F\colon\fromto{C}{\Top}$ just in case it is equivalent to the left fibration
\[\fromto{\Top_{\ast}\times_{\Top,F}C}{C}.\]
Dually, a right fibration $q\colon\fromto{Y}{C}$ is said to be \emph{classified by} $G\colon\fromto{C^{\op}}{\Top}$ just in case it is equivalent to the right fibration
\[\fromto{\Top^{\op}_{\ast}\times_{\Top^{\op},G}C}{C}.\]

The proofs of Joyal's theorem\footnote{We know four proofs: the original one due to Joyal, a modification thereof due to Lurie \cite{HTT}, a simplification due to Dugger and Spivak \cite{Dugger2009}, and a recent extreme simplification due to Heuts and Moerdijk \cite{HM1,HM2}.} are all relatively involved, and they involve breaking this assertion up into several constituent parts. But rather than get distracted by these details (beautiful though they be!), let us instead swim in the waters of appreciation for this result as end-users.

\begin{exm} Even for $C=\Delta^0$, this theorem is nontrivial: it provides an equivalence $\equivto{\Top}{\categ{Gpd}_{\infty}}$, where of course $\categ{Gpd}_{\infty}$ is the simplicial nerve of the full simplicial subcategory of $s\Set$ spanned by the Kan complexes. In other words, this result provides a concrete model for the $\infty$-category that was might have only been known through its universal property.
\end{exm}

But the deeper point here is that with the universal characterization of the introduction, it's completely unclear how to specify a functor from an $\infty$-category $C$ \emph{into} $\Top$. Even with the equivalence $\Top\simeq\categ{Gpd}_{\infty}$, we would still have to specify an infinite hierarchy of data to check this. However, with Joyal's result in hand, our task becomes to construct a left fibration $\fromto{X}{C}$. This is infinitely simpler: we are free to construct $X$ in a manner that will visibly map to our $C$, and our only job is to check a bunch of horn-filling conditions.

\begin{exm} One may use Joyal's theorem to find that if $p\colon\fromto{X}{C}$ is a left fibration classified by a functor $F\colon\fromto{C}{\Top}$, then the colimit of $F$ is weakly homotopy equivalent to $X$, and the the limit of $F$ is weakly homotopy equivalent to the space $\Map_{C}(C,X)$ of sections of $p$.
\end{exm}

\begin{exm} One attitude toward $\infty$-categories is that they are meant to be categories ``weakly enriched'' in spaces. Whatever this means, it should at least entail corepresentable and representable functors
\[h^x\colon\fromto{C}{\Top}\text{\quad and\quad}h_x\colon\fromto{C^{\op}}{\Top}.\]
But thanks to Joyal's theorem, we already have these: the former is given by the left fibration $\fromto{C_{x/}}{C}$, and the latter is given by the right fibration $\fromto{C_{/x}}{C}$.

This also provides a recognition principle: a left fibration $\fromto{X}{C}$ corresponds to a corepresentable functor if and only if $X$ admits an initial object; in this case, we call the left fibration itself \emph{corepresentable}. Dually, a right fibration $\fromto{X}{C}$ corresponds to a representable functor if and only if $X$ admits an terminal object; in this case, we call the right fibration itself \emph{representable}.
\end{exm}

\begin{exm} In the same vein, we expect to have a functor
\[\Map_C\colon\fromto{C^{\op}\times C}{\Top}\]
for any $\infty$-category $C$. From Joyal's Theorem, our job becomes to construct a left fibration $(s,t)\colon\fromto{\widetilde{\mathscr{O}}(C)}{C^{\op}\times C}$. It turns out that this isn't so difficult: define $\widetilde{\mathscr{O}}(C)$ via the formula
\[\widetilde{\mathscr{O}}(C)_n\coloneq\Mor(\Delta^{n,\op}\star\Delta^n,C),\]
and $s$ and $t$ are induced by the inclusions $\into{\Delta^{n,\op}}{\Delta^{n,\op}\star\Delta^n}$ and $\into{\Delta^n}{\Delta^{n,\op}\star\Delta^n}$, respectively.

So the claim is that $(s,t)$ is a left fibration. This isn't a completely trivial matter, but there is a proof in \cite{DAGX}, and another, slightly simpler, proof in \cite{Aefffib}. The key point is to study the behavior of the left adjoint of the functor $\widetilde{\mathscr{O}}$ on certain ``left anodyne'' monomorphisms.
\end{exm}

These two examples illustrate nicely a general principle about working ``vertically'' -- i.e., with left and right fibrations -- versus working ``horizontally'' -- i.e., with functors to $\Top$. It is easy to write down a left or right fibration, but it may not be easy to see that it \emph{is} a left or right fibration. On the other hand, it is quite difficult even to write down a suitable functor to $\Top$. So working vertically rather than horizontally relocates the difficulty in higher category theory from a struggle to make good definitions to a struggle to prove good properties.

\begin{exm} Suppose $C$ an $\infty$-category. A Kan fibration to $C$ is simultaneously a left fibration and a right fibration. So a Kan fibration must correspond to a both a covariant functor and a contravariant functor to $\Top$, and one sees that the ``pushforward'' maps must be homotopy inverse to the ``pullback'' maps. That is, the following are equivalent for a map $p\colon\fromto{X}{C}$ of simplicial sets:
\begin{itemize}
\item $p$ is a Kan fibration;
\item $p$ is a left fibration, and the functor $\fromto{C}{\Top}$ that classifies it carries any morphism of $C$ to an equivalence;
\item $p$ is a right fibration, and the functor $\fromto{C^{\op}}{\Top}$ that classifies it carries any morphism of $C$ to an equivalence.
\end{itemize}

From this point view, we see that when $C$ is an $\infty$-groupoid, Kan fibrations $\fromto{X}{C}$ are ``essentially the same thing'' as functors $\fromto{C}{\Top}$, which are in turn indistinguishable from functors $\fromto{C^{\op}}{\Top}$.
\end{exm}

\section{Inner fibrations and isofibrations} Inner fibrations are tricky to motivate from a $1$-categorical standpoint, because the nerve of any functor is automatically an inner fibration. We will discuss here a reasonable way of thinking about inner fibrations, but we do not know a reference for complete proofs, yet.

\begin{dfn} An \emph{inner fibration} is a map $p\colon\fromto{X}{S}$ of simplicial sets such that for any integer $n\geq 2$ and any $0<k<n$ and any solid arrow commutative square
\begin{equation*}
\begin{tikzpicture}[baseline]
\matrix(m)[matrix of math nodes,
row sep=4ex, column sep=4ex,
text height=1.5ex, text depth=0.25ex]
{\Lambda^n_k & X \\
\Delta^n & S, \\ };
\path[>=stealth,->,font=\scriptsize]
(m-1-1) edge node[above]{} (m-1-2)
edge[right hook->] (m-2-1)
(m-1-2) edge node[right]{$p$} (m-2-2)
(m-2-1) edge node[below]{} (m-2-2)
(m-2-1) edge[dotted] (m-1-2);
\end{tikzpicture}
\end{equation*}
there exists a dotted lift.
\end{dfn}

\begin{exm} Of course a simplicial set $X$ is an $\infty$-category just in case the canonical map $\fromto{X}{\Delta^0}$ is an inner fibration. Consequently, any fiber of an inner fibration is an $\infty$-category.
\end{exm}

\begin{exm} If $X$ is an $\infty$-category and $D$ is an ordinary category, then it's easy to see that any map $\fromto{X}{ND}$ is an inner fibration.

On the other hand, a map $p\colon\fromto{X}{S}$ is an inner fibration if and only if, for any $n$-simplex $\sigma\in S_n$, the pullback
\[\fromto{X\times_{S,\sigma}\Delta^n}{\Delta^n}\]
is an inner fibration. Consequently, we see that $p$ is an inner fibration if and only if, for any $n$-simplex $\sigma\in S_n$, the pullback $X\times_{S,\sigma}\Delta^n$ is an $\infty$-category.
\end{exm}

So in a strong sense, we'll understand the ``meaning'' of inner fibrations one we understand the ``meaning'' of functors from $\infty$-categories to $\Delta^n$.

\begin{exm} When $n=1$, we have the following. For any $\infty$-categories $C_0$ and $C_1$, there is an equivalence of $\infty$-categories
\[\equivto{\{C_0\}\times_{\Cat_{\infty/\Delta^{\{0\}}}}\Cat_{\infty/\Delta^1}\times_{\Cat_{\infty/\Delta^{\{1\}}}}\{C_1\}}{\Fun(C_0^{\op}\times C_1,\Top)}.\]
A proof of this fact doesn't seem to be contained in the literature yet, but we will nevertheless take it as given.

Now the $\infty$-category on the right of this equivalence can also be identified with the $\infty$-category
\[\Fun^L(P(C_1),P(C_0))\]
of colimit-preserving functors between $P(C_1)=\Fun(C_1^{\op},\Top)$ and $P(C_0)=\Fun(C_0^{\op},\Top)$. Such a colimit-preserving functor is sometimes called a \emph{profunctor}.
\end{exm}

\begin{exm} When $n=2$, if $C$ is an $\infty$-category, and $\fromto{C}{\Delta^2}$ is a functor, then we have three fibers $C_0$, $C_1$, and $C_2$ and three colimit-preserving functors $F\colon\fromto{P(C_2)}{P(C_1)}$, $G\colon\fromto{P(C_1)}{P(C_0)}$, and $H\colon\fromto{P(C_2)}{P(C_0)}$. Furthermore, there is natural transformation $\alpha\colon\fromto{G\circ F}{H}$.
\end{exm}

In the general case, a functor $\fromto{C}{\Delta^n}$ amounts to the choice of $\infty$-categories $C_0, C_1,\dots, C_n$ and a lax-commutative diagram of colimit-preserving functors among the various $\infty$-cate\-gories $P(C_i)$.

What we would like to say now is that the $\infty$-category of inner fibrations $\fromto{X}{S}$ is equivalent to the $\infty$-category of lax functors from $S^{\op}$ to a suitable ``double $\infty$-category'' of $\infty$-categories and profunctors. We do not know, however, how to make such an assertion precise.

In any case, we could ask for more restrictive hypotheses. We could, for example, ask for fibrations $\fromto{X}{S}$ that are classified by functors from $S$ to an $\infty$-category of profunctors (so that all the $2$-morphisms that appear are equivalences); this is covered by the notion of \emph{flatness}, which will discuss in the section after next. More restrictively, we can ask for fibrations $\fromto{X}{S}$ that are classified by functors from $S$ to $\Cat_{\infty}$ itself; these are \emph{cocartesian fibrations}, which we will discuss in the next section.

For future reference, let's specify an extremely well-behaved class of inner fibrations.
\begin{dfn} Suppose $C$ an $\infty$-category. Then an \emph{isofibration} (\textsc{aka} a \emph{categorical fibration}\footnote{Emily Riehl makes the clearly compelling case that ``isofibration'' is preferable terminology, because it actually suggests what kind of lifting property it will have, whereas the word ``categorical'' is unhelpful in this regard. She also tells us that ``isofibration'' is a standard term in $1$-category theory, and that the nerve of functor is an isofibration iff the functor is an isofibration. We join her in her view that ``isofibration'' is better.}) $p\colon\fromto{X}{C}$ is an inner fibration such that for any object $x\in X_0$ and any equivalence $f\colon\fromto{p(x)}{b}$ of $C$, there exists an equivalence $\phi\colon\fromto{x}{y}$ of $X$ such that $p(\phi)=f$.
\end{dfn}

\noindent We shall revisit this notion in greater detail in a moment, but for now, let us simply comment that an isofibration $\fromto{X}{C}$ is an inner fibration whose fibers vary functorially in the equivalences of $C$.

For more general bases, this definition won't do, of course, but we won't have any use for isofibrations whose targets are not $\infty$-categories.

\section{Cocartesian and cartesian fibrations} If $F\colon\fromto{C}{\Cat}$ is an (honest) diagram of ordinary categories, then one can generalize the category of elements construction as follows: form the category $X$ whose objects are pairs $(c,x)$ consisting of an object $c\in C$ and an object $x\in F(c)$, in which a morphism $(f,\phi)\colon\fromto{(d,y)}{(c,x)}$ is a morphism $f\colon\fromto{d}{c}$ of $C$ and a morphism
\begin{equation*}
\phi\colon\fromto{X(f)(y)}{x}
\end{equation*}
of $F(c)$. This is called the \emph{Grothendieck construction}, and there is an obvious forgetful functor $p\colon\fromto{X}{C}$.

One may now attempt to reverse-engineer the Grothendieck construction by trying to extract the salient features of the forgetful functor $p$. What we may notice is that for any morphism $f\colon\fromto{d}{c}$ of $C$ and any object $y\in F(d)$ there is a special morphism
\begin{equation*}
\Phi=(f,\phi)\colon\fromto{(d,y)}{(c,X(f)(y))}
\end{equation*}
of $X$ in which
\begin{equation*}
\phi\colon\fromto{F(f)(y)}{F(f)(y)}
\end{equation*}
is simply the identity morphism. This morphism is \emph{initial} among all the morphisms $\Psi$ of $X$ such that $p(\Psi)=f$; that is, for any morphism $\Psi$ of $X$ such that $p(\Psi)=f$, there exists a morphism $I$ of $X$ such that $p(I)=\id_c$ such that $\Psi=I\circ\Phi$.

We call morphisms of $X$ that are initial in this sense \emph{$p$-cocartesian}. Since a $p$-cocartesian edge lying over a morphism $\fromto{d}{c}$ is defined by a universal property, it is uniquely specified up to a unique isomorphism lying over $\id_{c}$. The key condition that we are looking for is then that \emph{for any morphism of $C$ and any lift of its source, there is a $p$-cocartesian morphism with that source lying over it}. A functor $p$ satisfying this condition is called a \emph{Grothendieck opfibration}.

Now for any Grothendieck opfibration $p\colon\fromto{X}{C}$, let us attempt to extract a functor $F\colon\fromto{C}{\Cat}$ that gives rise to it in this way. We proceed in the following manner. To any object $a\in C$ assign the fiber $X_a$ of $p$ over $a$. To any morphism $f\colon\fromto{a}{b}$ assign a functor $F(f)\colon\fromto{X_a}{X_b}$ that carries any object $x\in X_a$ to the target $F(f)(x)\in X_b$ of ``the'' $q$-cocartesian edge lying over $f$.

Right away, we have a problem: $q$-cocartesian edges are only unique up to isomorphism. So these functors cannot be strictly compatible with composition; rather, one will obtain natural isomorphisms
\begin{equation*}
F(g\circ f)\simeq F(g)\circ F(f)
\end{equation*}
that will satisfy a secondary layer of coherences that make $F$ into a \emph{pseudofunctor}. Fortunately, one can rectify this pseudofunctor to an equivalent honest functor, which in turn gives rise to $p$, up to equivalence.

As we have seen in our discussion of left and right fibrations, there are genuine advantages in homotopy theory to working with fibrations instead of functors. Consequently, we define a class of fibrations that is a natural generalization of the class of Grothendieck opfibrations.
\begin{dfn}\label{rec:cocart} If $p\colon\fromto{X}{S}$ is an inner fibration of simplicial sets, then an edge $f\colon\fromto{\Delta^{1}}{X}$ is \emph{$p$-cocartesian} just in case, for each integer $n\geq 2$, any extension
\begin{equation*}
\begin{tikzpicture} 
\matrix(m)[matrix of math nodes, 
row sep=4ex, column sep=4ex, 
text height=1.5ex, text depth=0.25ex] 
{\Delta^{\{0,1\}}&X,\\ 
\Lambda^n_0&\\}; 
\path[>=stealth,->,font=\scriptsize] 
(m-1-1) edge node[above]{$f$} (m-1-2) 
edge[right hook->] (m-2-1) 
(m-2-1) edge node[below]{$F$} (m-1-2); 
\end{tikzpicture}
\end{equation*}
and any solid arrow commutative diagram
\begin{equation*}
\begin{tikzpicture} 
\matrix(m)[matrix of math nodes, 
row sep=4ex, column sep=4ex, 
text height=1.5ex, text depth=0.25ex] 
{\Lambda^n_0&X\\ 
\Delta^n&S,\\}; 
\path[>=stealth,->,font=\scriptsize] 
(m-1-1) edge node[above]{$F$} (m-1-2) 
edge[right hook->] (m-2-1) 
(m-1-2) edge node[right]{$p$} (m-2-2) 
(m-2-1) edge (m-2-2)
(m-2-1) edge[dotted] node[below]{$\overline{F}$} (m-1-2);
\end{tikzpicture}
\end{equation*}
a dotted lift exists.

We say that $p$ is a \emph{cocartesian fibration} if, for any edge $\eta\colon\fromto{s}{t}$ of $S$ and for every vertex $x\in X_0$ such that $p(x)=s$, there exists a $p$-cocartesian edge $f\colon\fromto{x}{y}$ such that $\eta=p(f)$.

\emph{Cartesian edges} and \emph{cartesian fibrations} are defined dually, so that an edge of $X$ is $p$-cartesian just in case the corresponding edge of $X^{\op}$ is cocartesian for the inner fibration $p^{\op}\colon\fromto{X^{\op}}{S^{\op}}$, and $p$ is a cartesian fibration just in case $p^{\op}$ is a cocartesian fibration.
\end{dfn}

\begin{exm}[\protect{\cite[Rk 2.4.2.2]{HTT}}]\label{exm:nerveofGrothfib} A functor $p\colon\fromto{D}{C}$ between ordinary categories is a Grothendieck opfibration if and only if the induced functor $N(p)\colon\fromto{ND}{NC}$ on nerves is a cocartesian fibration.
\end{exm}

\begin{exm} Any left fibration is a cocartesian fibration, and a cocartesian fibration is a left fibration just in case its fibers are $\infty$-groupoids.

Dually, of course, the class of right fibrations coincides with the class of cartesian fibrations whose fibers are $\infty$-groupoids.
\end{exm}

\begin{exm} Suppose $C$ an $\infty$-category and $p\colon\fromto{X}{C}$ an inner fibration. Then for any morphism $\eta$ of $X$, the following are equivalent.
\begin{itemize}
\item $\eta$ is an equivalence of $X$;
\item $\eta$ is $p$-cocartesian, and $p(\eta)$ is an equivalence of $C$;
\item $\eta$ is $p$-cartesian, and $p(\eta)$ is an equivalence of $C$.
\end{itemize}
It follows readily that if $p$ is a cocartesian or cartesian fibration, then it is an isofibration.

Conversely, if $C$ is an $\infty$-groupoid, then the following are equivalent.
\begin{itemize}
\item $p$ is an isofibration;
\item $p$ is a cocartesian fibration;
\item $p$ is a cartesian fibration.
\end{itemize}
\end{exm}

\begin{exm}[\protect{\cite[Cor. 2.4.7.12]{HTT}}]\label{exm:scocartfib} For any $\infty$-category $C$, we write $\mathscr{O}(C)\coloneq\Fun(\Delta^1,C)$. Evaluation at $0$ defines a cartesian fibration $s\colon\fromto{\mathscr{O}(C)}{C}$, and evaluation at $1$ defines a cocartesian fibration $t\colon\fromto{\mathscr{O}(C)}{C}$.

One can ask whether the functor $s\colon\fromto{\mathscr{O}(C)}{C}$ is also a \emph{cocartesian} fibration. One may observe \cite[Lm. 6.1.1.1]{HTT} that an edge $\fromto{\Delta^1}{\mathscr{O}(C)}$ is $s$-cocartesian just in case the correponding diagram $\fromto{(\Lambda^2_0)^{\rhd}\cong\Delta^1\times\Delta^1}{C}$ is a pushout square.
\end{exm}

\begin{exm} Consider the full subcategory $\categ{RFib}^{\textit{rep}}\subset\mathscr{O}(\Cat_{\infty})$ spanned by the representable right fibrations. The restriction of the functor $t\colon\fromto{\mathscr{O}(\Cat_{\infty})}{\Cat_{\infty}}$ to $\categ{RFib}^{\textit{rep}}$ is again a cocartesian fibration.
\end{exm}

In the following, we will denote by $\Cat_{\infty}$ the simplicial nerve of the (fibrant) simplicial category whose objects are $\infty$-categories, in which $\Map(C,D)$ is the maximal $\infty$-groupoid contained in $\Fun(C,D)$. Similarly, for any $\infty$-category $C$, we will denote by $\categ{Cocart}(C)$ (respectively, $\categ{Cart}(C)$) the simplicial nerve of the (fibrant) simplicial category whose objects are cocartesian (resp., cartesian) fibrations $\fromto{X}{C}$, in which $\Map(X,Y)$ is the $\infty$-groupoid whose $n$-simplices are functors $\fromto{X\times\Delta^n}{Y}$ over $C$ that carry any edge $(f,\tau)$ in which $f$ is cocartesian (resp. cartesian) to a cocartesian edge (resp., a cartesian edge).

\begin{thm} Suppose $C$ an $\infty$-category. For any functor $F\colon\fromto{C}{\Cat_{\infty}}$, we may consider the cocartesian fibration
\[\fromto{\categ{RFib}^{\textit{rep}}\times_{\Cat_{\infty},F}C}{C}.\]
This defines an equivalence of categories
\[\equivto{\Fun(C,\Cat_{\infty})}{\categ{Cocart}(C)}.\]

Dually, for any functor $G\colon\fromto{C}{\Cat_{\infty}^{\op}}$, we may consider the cartesian fibration
\[\fromto{\categ{RFib}^{\textit{rep},\op}\times_{\Cat_{\infty}^{\op},G}C}{C}.\]
This defines an equivalence of categories
\[\equivto{\Fun(C^{\op},\Cat_{\infty})}{\categ{Cart}(C)}.\]
\end{thm}

\noindent A cocartesian fibration $p\colon\fromto{X}{C}$ is said to be \emph{classified} by $F$ just in case it is equivalent to the cocartesian fibration
\[\fromto{\categ{RFib}^{\textit{rep}}\times_{\Cat_{\infty},F}C}{C}.\]
Dually, a cartesian fibration $q\colon\fromto{Y}{C}$ is said to be \emph{classified} by $F$ just in case it is equivalent to the cocartesian fibration
\[\fromto{\categ{RFib}^{\textit{rep},\op}\times_{\Cat_{\infty}^{\op},G}C}{C}.\]

\begin{exm} Suppose $C$ an $\infty$-category, and suppose $\fromto{X}{C}$ an isofibration. If $\iota C\subseteq C$ is the largest $\infty$-groupoid contained in $C$, then the pulled back isofibration
\[\fromto{X\times_C\iota C}{\iota C}\]
is both cocartesian and cartesian, and so it corresponds to a functor
\[\fromto{\iota C\simeq\iota C^{\op}}{\Cat_{\infty}}.\]
This is the sense in which the fibers of an isofibration vary functorially in equivalences if $C$.
\end{exm}

\begin{exm} For any $\infty$-category $C$, the functor $\fromto{C^{\op}}{\Cat_{\infty}}$ that classifies the cartesian fibration $s\colon\fromto{\mathscr{O}(C)}{C}$ is the functor that carries any object $a$ of $C$ to the undercategory $C_{a/}$ and any morphism $f\colon\fromto{a}{b}$ to the forgetful functor $f^{\star}\colon\fromto{C_{b/}}{C_{a/}}$.

If $C$ admits all pushouts, then the cocartesian fibration $s\colon\fromto{\mathscr{O}(C)}{C}$ is classified by a functor $\fromto{C}{\Cat_{\infty}}$ that carries any object $a$ of $C$ to the undercategory $C_{a/}$ and any morphism $f\colon\fromto{a}{b}$ to the functor $f_{!}\colon\fromto{C_{a/}}{C_{b/}}$ that is given by pushout along $f$.
\end{exm}

One particularly powerful construction with cartesian and cocartesian fibrations comes from \cite[\S 3.2.2]{HTT}. We've come to call this the \emph{cartesian workhorse.}
\begin{exm}\label{cor32213} Suppose $p\colon\fromto{X}{B^{\op}}$ a cartesian fibration and $q\colon\fromto{Y}{B^{\op}}$ a cocartesian fibration. Suppose $F\colon\fromto{B}{\Cat_{\infty}}$ a functor that classifies $p$ and $G\colon\fromto{B^{\op}}{\Cat_{\infty}}$ a functor that classifies $q$. Clearly one may define a functor
\begin{equation*}
\Fun(F,G)\colon\fromto{B^{\op}}{\Cat_{\infty}}
\end{equation*}
that carries a vertex $s$ of $B^{\op}$ to the $\infty$-category $\Fun(F(s),G(s))$ and an edge $\eta\colon\fromto{s}{t}$ of $B^{\op}$ to the functor
\begin{equation*}
\fromto{\Fun(F(s),G(s))}{\Fun(F(t),G(t))}
\end{equation*}
given by the assignment $\goesto{F}{G(\eta)\circ F\circ F(\eta)}$.

If one wishes to work instead with the fibrations directly (avoiding straightening and unstraightening), the following construction provides an elegant way of writing explicitly the cocartesian fibration classified by the functor $\Fun(F,G)$.

Suppose $p\colon\fromto{X}{B^{\op}}$ is a cartesian fibration classified by a functor
\[F\colon\fromto{B}{\Cat_{\infty}},\]
and suppose $q\colon\fromto{Y}{B^{\op}}$ is a cocartesian fibration classified by a functor
\[G\colon\fromto{B^{\op}}{\Cat_{\infty}}.\]
One defines a simplicial set $\widetilde{\underline{\Fun}}_B(X,Y)$ and a map $r\colon\fromto{\widetilde{\underline{\Fun}}_B(X,Y)}{B^{\op}}$ defined by the following universal property: for any map $\sigma\colon\fromto{K}{B^{\op}}$, one has a bijection
\begin{equation*}
\Mor_{/B^{\op}}(K,\widetilde{\underline{\Fun}}_B(X,Y))\cong\Mor_{/B^{\op}}(X\times_{B^{\op}}K,Y),
\end{equation*}
functorial in $\sigma$.

It is then shown in \cite[Cor. 3.2.2.13]{HTT} (but see \ref{exm:cartesianworkhorse} below for a better proof) that $r$ is a cocartesian fibration, and an edge
\[g\colon\fromto{\Delta^1}{\widetilde{\underline{\Fun}}_B(X,Y)}\]
is $r$-cocartesian just in case the induced map $\fromto{X\times_{B^{\op}}\Delta^1}{Y}$ carries $p$-cartesian edges to $q$-cocartesian edges. The fiber of the map $\fromto{\widetilde{\underline{\Fun}}_B(X,Y)}{S}$ over a vertex $s$ is the $\infty$-category $\Fun(X_s,Y_s)$, and for any edge $\eta\colon\fromto{s}{t}$ of $B^{\op}$, the functor $\eta_!\colon\fromto{T_s}{T_t}$ induced by $\eta$ is equivalent to the functor $\goesto{F}{G(\eta)\circ F\circ F(\eta)}$ described above.
\end{exm}

\section{Flat inner fibrations} We have observed that if $C$ is an $\infty$-category, and if $\fromto{C}{\Delta^1}$ is any functor with fibers $C_0$ and $C_1$, then there is a corresponding profunctor from $C_1$ to $C_0$, i.e., a colimit-preserving functor $\fromto{P(C_1)}{P(C_0)}$. Furthermore, the passage from $\infty$-categories over $\Delta^1$ to profunctors is even in some sense an equivalence.

So to make this precise, let $\categ{Prof}$ denote the full subcategory of the $\infty$-category $\categ{Pr}^L$ of presentable $\infty$-categories and left adjoints spanned by those $\infty$-categories of the form $P(C)$. We can \emph{almost} -- but not quite -- construct an equivalence between the $\infty$-categories $\Cat_{\infty/\Delta^1}$ and the $\infty$-category $\Fun(\Delta^1,\categ{Prof})$.

The trouble here is that there are strictly more equivalences in $\categ{Prof}$ than there are in $\Cat_{\infty}$; two $\infty$-categories are equivalent in $\categ{Prof}$ if and only if they have equivalent idempotent completions. So that suggests the fix for this problem: we employ a pullback that will retain the data of the $\infty$-categories that are the source and target of our profunctor.

\begin{prp} There is an equivalence of $\infty$-categories
\[\Cat_{\infty/\Delta^1}\simeq\Fun(\Delta^1,\categ{Prof})\times_{\Fun(\partial\Delta^1,\categ{Prof})}\Fun(\partial\Delta^1,\Cat_{\infty}).\]
\end{prp}
\noindent We do not know of a reference for this result, yet, but we expect this to appear in a future work of P. Haine.

When we pass to $\infty$-categories over $\Delta^2$, we have a more complicated problem: a functor $\fromto{C}{\Delta^2}$ only specifies a lax commutative diagram of profunctors: three fibers $C_0$, $C_1$, and $C_2$; three colimit-preserving functors $F\colon\fromto{P(C_0)}{P(C_1)}$, $G\colon\fromto{P(C_1)}{P(C_2)}$, and $H\colon\fromto{P(C_0)}{P(C_2)}$; and a natural transformation $\alpha\colon\fromto{G\circ F}{H}$. In order to ensure that $\alpha$ be a natural equivalence, we need a condition on our fibration. This is where flatness comes in.

\begin{dfn}\label{dfn:flatinnerfib} An inner fibration $p\colon\fromto{X}{S}$ is said to be \emph{flat} just in case, for any inner anodyne map $\into{K}{L}$ and any map $\fromto{L}{S}$, the pullback
\[\fromto{X\times_SK}{X\times_SL}\]
is a categorical equivalence.
\end{dfn}

We will focus mostly on flat \emph{categorical} fibrations. Let us see right away that some familiar examples and constructions yield flat inner fibrations.

\begin{exm}[Lurie, \protect{\cite[Ex. B.3.11]{HA}}] Cocartesian and cartesian fibrations are flat isofibrations. In particular, if $C$ is an $\infty$-groupoid, then any isofibration $\fromto{X}{C}$ is flat.
\end{exm}

\begin{exm}[Lurie, \protect{\cite[Pr. B.3.13]{HA}}] If $p\colon\fromto{X}{S}$ is a flat inner fibration, then for any vertex $x\in X_0$, the inner fibrations $\fromto{X_{x/}}{S_{p(x)/}}$ and $\fromto{X_{/x}}{S_{/p(x)}}$ are flat as well.
\end{exm}

It is not necessary to test flatness with \emph{all} inner anodyne maps; in fact, one can make do with the inner horn of a $2$-simplex:
\begin{prp}[Lurie, \protect{\cite[Pr. B.3.14]{HA}}] An inner fibration $p\colon\fromto{X}{S}$ is flat if and only if, for any $2$-simplex $\sigma\in S_2$, the pullback
\[\fromto{X\times_{p,S,(\sigma|\Lambda^2_1)}\Lambda^2_1}{X\times_{p,S,\sigma}\Delta^2}\]
is a categorical equivalence.
\end{prp}

\begin{prp}[Lurie, \protect{\cite[Pr. B.3.2, Rk. B.3.9]{HA}}] An inner fibration $\fromto{X}{S}$ is flat just in case, for any $2$-simplex
\begin{equation*}
\begin{tikzpicture}[baseline]
\matrix(m)[matrix of math nodes,
row sep=4ex, column sep=4ex,
text height=1.5ex, text depth=0.25ex]
{&v&\\
u&&w\\};
\path[>=stealth,->,font=\scriptsize]
(m-1-2) edge (m-2-3)
(m-2-1) edge (m-1-2)
edge node[below]{$f$} (m-2-3);
\end{tikzpicture}
\end{equation*}
of $S$ any for any edge $\fromto{x}{y}$ lying over $f$, the simplicial set $X_{x/\ /y}\times_S\{v\}$ is weakly contractible.
\end{prp}

\begin{nul} Note that the condition of the previous result is vacuous if the $2$-simplex is degenerate. Consequently, if $S$ is $1$-skeletal, then any inner fibration $\fromto{X}{S}$ is flat.
\end{nul}

\begin{prp} Suppose $C$ an $\infty$-category. Then there is an equivalence
\[\equivto{\categ{Flat}(C)}{\Fun(C^{\op},\categ{Prof})\times_{\Fun(\iota C^{\op},\categ{Prof})}\Fun(\iota C^{\op},\Cat_{\infty})},\]
where $\categ{Flat}(C)$ is an $\infty$-category of flat isofibrations $\fromto{X}{C}$.
\end{prp}

\noindent Once again, we do not know a reference for this in the literature yet, but we expect that this will be shown in future work of P. Haine.

\section{Categorical patterns} One of the primary appeals that flat inner fibrations have to offer is a collection of right Quillen functors for certain combinatorial model structures defined by means of \emph{categorical patterns}.

\begin{dfn} A \emph{categorical pattern} on a simplicial set $S$ is a triple $(M,T,P)$ consisting of the following:
\begin{itemize}
\item a set $M\subset S_1$ of \emph{marked edges} that contains all the degenerate edges,
\item a set $T\subset S_2$ of \emph{scaled $2$-simplices} that contains all the degenerate $2$-simplices, and
\item a set $P$ of maps $f_{\alpha}\colon\fromto{K_{\alpha}^{\lhd}}{S}$ such that $f((K_{\alpha}^{\lhd})_1)\subset M$ and $f((K_{\alpha}^{\lhd})_2)\subset T$.
\end{itemize}
\end{dfn}

\begin{ntn} For the purposes of this appendix, if $p\colon\fromto{X}{S}$ and $f\colon\fromto{K}{S}$ are maps of simplicial sets, then let us write
\[p_f\colon\fromto{X\times_SK}{K}\]
for the pullback of $p$ along $f$.
\end{ntn}

\begin{dfn} Suppose $(M,T,P)$ a categorical pattern on a simplicial set $S$. Then a marked map $p\colon\fromto{(X,E)}{(S,M)}$ is said to be \emph{$(M,T,P)$-fibered} if the following conditions obtain.
\begin{itemize}
\item The map $p\colon\fromto{X}{S}$ is an inner fibration.
\item For every marked edge $\eta$ of $S$, the pullback $p_{\eta}$ is a cocartesian fibration.
\item An edge $\epsilon$ of $X$ is marked just in case it is $p_{p(\epsilon)}$-cocartesian.
\item For any commutative square
\begin{equation*}
\begin{tikzpicture}[baseline]
\matrix(m)[matrix of math nodes,
row sep=4ex, column sep=4ex,
text height=1.5ex, text depth=0.25ex]
{\Delta^{\{0,1\}} & X\\
\Delta^2 &  S\\ };
\path[>=stealth,->,font=\scriptsize]
(m-1-1) edge node[above]{$\epsilon$} (m-1-2)
edge[right hook->] (m-2-1)
(m-1-2) edge node[right]{$p$} (m-2-2)
(m-2-1) edge node[below]{$\sigma$} (m-2-2);
\end{tikzpicture}
\end{equation*}
in which $\sigma$ is scaled, if $\epsilon$ is marked, then it is $p_{\sigma}$-cocartesian.
\item For every element $f_{\alpha}\colon\fromto{K_{\alpha}^{\lhd}}{S}$ of $P$, the cocartesian fibration
\[p_{f_{\alpha}}\colon\fromto{X\times_SK_{\alpha}^{\lhd}}{K_{\alpha}^{\lhd}}\]
is classified by a limit diagram $\fromto{K_{\alpha}^{\lhd}}{\Cat_{\infty}}$.
\item For every element $f_{\alpha}\colon\fromto{K_{\alpha}^{\lhd}}{S}$ of $P$, any cocartesian section of $p_{f_{\alpha}}$ is a $p$-limit diagram in $X$.
\end{itemize}
\end{dfn}

\begin{exm} For any simplicial set $S$, the $(S_1,S_2,\varnothing)$-fibered maps $\fromto{(X,E)}{S^{\sharp}}$ are precisely those maps of the form $\fromto{\leftnat{X}}{S^{\sharp}}$, where the underlying map of simplicial sets $\fromto{X}{S}$ is a cocartesian fibration.

In particular, if $B$ is an orbital $\infty$-category, then an $(S_1,S_2,\varnothing)$-fibered maps $\fromto{(X,E)}{S^{\sharp}}$ is essentially the same thing as a $B$-$\infty$-category.
\end{exm}

\begin{exm} Suppose $\Phi$ a perfect operator category \cite[\S 6]{opcat}, and consider the following categorical pattern
\begin{equation*}
(\mathrm{Ne},N(\Lambda(\Phi))_2,P)
\end{equation*}
on the nerve $N\Lambda(\Phi)$ of the Leinster category \cite[\S 7]{opcat}. Here the class $\mathrm{Ne}\subset N(\Lambda(\Phi))_1$ consists of all the inert morphisms of $N\Lambda(\Phi)$. The class $P$ is the set of maps $\fromto{\Lambda^2_0}{N\Lambda(\Phi)}$ given by diagrams $I\ot J\to I'$ of $\Lambda(\Phi)$ in which both $J\to I$ and $J\to I'$ are inert, and
\begin{equation*}
|J|=|J\times_{TI}I|\sqcup|J\times_{TI'}I'|\cong|I|\sqcup|I'|.
\end{equation*}

Then a marked map $\fromto{(X,E)}{(N\Lambda(\Phi),N\Lambda^{\dag}(\Phi))}$ is $(\mathrm{Ne},N(\Lambda(\Phi))_2,P)$-fibered just in case the underlying map of simplicial sets $\fromto{X}{N\Lambda(\Phi)}$ is a $\infty$-operad over $\Phi$, and $E$ is the collection of cocartesian edges over the inert edges \cite[\S 8]{opcat}.

In particular, when $\Phi=\FF$, the category of finite sets, this recovers the collection of $\infty$-operads; cf. \cite[Pr. 2.4.1.6]{HA}.
\end{exm}

\begin{exm} If $\Phi$ is a perfect operator category, we can contemplate another categorical pattern
\begin{equation*}
(N(\Lambda(\Phi))_1,N(\Lambda(\Phi))_2,P)
\end{equation*}
on $\Lambda(\Phi)$. Here $P$ is as in the previous example.

Now a marked map $\fromto{(X,E)}{(N\Lambda(\Phi),N\Lambda^{\dag}(\Phi))}$ is $(N(\Lambda(\Phi))_1,N(\Lambda(\Phi))_2,P)$-fibered just in case the underlying map of simplicial sets $\fromto{X}{N\Lambda(\Phi)}$ is a $\Phi$-monoidal $\infty$-category, and $E$ is the collection of all cocartesian edges.
\end{exm}

\begin{thm}[Lurie, \protect{\cite[Th. B.0.20]{HA}}]\label{thm:modstructcatpatter} Suppose $(M,T,P)$ a categorical pattern on a simplicial set $S$. Then there exists a left proper, combinatorial, simplicial model structure on $s\Set^+_{/(S,M)}$ in which the cofibrations are monomorphisms, and a marked map $p\colon\fromto{(X,E)}{(S,M)}$ is fibrant just in case it is $(M,T,P)$-fibered.

We will denote this model category $s\Set^+_{/(S,M,T,P)}$.
\end{thm}

We also have a characterization of the fibrations between fibrant objects in this model structure:
\begin{prp}[Lurie, \protect{\cite[Th. B.2.7]{HA}}] Suppose $S$ an $\infty$-category, and suppose $(M,T,P)$ a categorical pattern on $S$ such that every equivalence of $S$ is marked, and every $2$-simplex $\fromto{\Delta^2}{S}$ whose restriction to $\Delta^{\{0,1\}}$ is an equivalence is scaled. If $(Y,F)$ is fibrant in $s\Set^+_{/(S,M,T,P)}$, then a map $f\colon\fromto{(X,E)}{(Y,F)}$ over $(S,M)$ is a fibration in $s\Set^+_{/(S,M,T,P)}$ if and only if $(X,E)$ is fibrant, and the underlying map $f\colon\fromto{X}{Y}$ is an isofibration.
\end{prp}

\begin{exm} If $(S,M)$ is a marked simplicial set, then a typical categorical pattern on a simplicial set $S$ is simply $(M,S_2,\varnothing)$. We call the resulting model structure on $s\Set^+_{/(S,M)}$ the \emph{cocartesian model structure}.

Indeed, the cocartesian model structure on $s\Set^+_{/S^{\sharp}}$ is precisely the cocartesian model structure described by Lurie. Furthermore, in light of the previous proposition, if $p\colon\fromto{C}{S}$ is a cocartesian fibration, then we claim that the cocartesian model structure on $s\Set^{+}_{/\leftnat{C}}$ is created by the forgetful functor to the cocartesian model structure $s\Set^+_{/S}$ -- i.e., by regarding $s\Set^+_{/\leftnat{C}}$ as $(s\Set^+_{/S})_{/p}$. (Note that a model structure is uniquely charcterized by its cofibrations and fibrant objects.)
\end{exm}

\section{Constructing fibrations via additional functoriality} Now we can use the theory of flat isofibrations to describe some extra functoriality of these categorical pattern model structures in $S$. We will use this extra direction of functoriality to perform some beautifully explicit constructions of fibrations of various kinds. The main result is a tad involved, but the constructive power it offers is worth it in the end.

\begin{ntn} Let us begin with the observation that, given a map of marked simplicial sets $\pi\colon\fromto{(S,M)}{(T,N)}$, there is a string of adjoints
\[\pi_!\dashv\pi^{\star}\dashv\pi_{\star},\]
where the far left adjoint
\[\pi_!\colon\fromto{s\Set^+_{/(S,M)}}{s\Set^+_{/(T,N)}}\]
is simply composition with $\pi$; the middle functor
\[\pi^{\star}\colon\fromto{s\Set^+_{/(T,N)}}{s\Set^+_{/(S,M)}}\]
is given by the assignment $\goesto{(Y,F)}{(Y,F)\times_{(T,N)}(S,M)}$; and the far right adjoint
\[\pi_{\star}\colon\fromto{s\Set^+_{/(S,M)}}{s\Set^+_{/(T,N)}}\]
is given by a ``space of sections.'' That is, an $n$-simplex of $\pi_{\star}(X,E)$ is a pair $(\sigma,f)$ consisting of an $n$-simplex $\sigma\in T_n$ along with a marked map
\[\fromto{(\Delta^n)^{\flat}\times_{(T,N)}(S,M)}{(X,E)}\]
over $(S,M)$; an edge $(\eta,f)$ of $\pi_{\star}(X,E)$ is marked just in case $\eta$ is marked, and $f$ carries any edge of $\Delta^1\times_TS$ that projects to a marked edge of $S$ to a marked edge of $X$.
\end{ntn}

\begin{thm}[Lurie, \protect{\cite[Th. B.4.2]{HA}}]\label{thm:functorialityofcatpatterns} Suppose $S$, $S'$, and $X$ three $\infty$-categories; suppose $(M,T,P)$ a categorical pattern on $S$; suppose $(M',T',P')$ a categorical pattern on $S'$; and suppose $E$ a collection of marked edges on $X$. Assume these data satisfy the following conditions.
\begin{itemize}
\item Any equivalence of either $S$ or $X$ is marked.
\item The marked edges in $S$ and $X$ are each closed under composition.
\item Every $2$-simplex $\fromto{\Delta^2}{S}$ whose restriction to $\Delta^{\{0,1\}}$ is an equivalence is scaled.
\end{itemize}
Suppose
\[\pi\colon\fromto{(X,E)}{(S,M)}\]
a map of marked simplicial sets satisfying the following conditions.
\begin{itemize}
\item The functor $\pi\colon\fromto{X}{S}$ is a flat isofibration.
\item For any marked edge $\eta$ of $S$, the pullback
\[\pi_{\eta}\colon\fromto{X\times_S\Delta^1}{\Delta^1}\]
is a cartesian fibration.
\item For any element $f_{\alpha}\colon\fromto{K_{\alpha}^{\lhd}}{S}$ of $P$, the simplicial set $K_{\alpha}$ is an $\infty$-category, and the pullback
\[\pi_{f_{\alpha}}\colon\fromto{X\times_SK_{\alpha}^{\lhd}}{K_{\alpha}^{\lhd}}\]
is a cocartesian fibration.
\item Suppose
\begin{equation*}
\begin{tikzpicture}[baseline]
\matrix(m)[matrix of math nodes,
row sep=4ex, column sep=4ex,
text height=1.5ex, text depth=0.25ex]
{&y&\\
x&&z\\};
\path[>=stealth,->,font=\scriptsize]
(m-1-2) edge[inner sep=0.75pt] node[above right]{$\psi$} (m-2-3)
(m-2-1) edge[inner sep=0.75pt] node[above left]{$\phi$} (m-1-2)
edge node[below]{$\chi$} (m-2-3);
\end{tikzpicture}
\end{equation*}
a $2$-simplex of $X$ in which $\pi(\phi)$ is an equivalence, $\psi$ is locally $\pi$-cartesian, and $\pi(\psi)$ is marked. Then the edge $\phi$ is marked just in case $\chi$ is.
\item Suppose $f_{\alpha}\colon\fromto{K_{\alpha}^{\lhd}}{S}$ an element of $P$, and suppose
\begin{equation*}
\begin{tikzpicture}[baseline]
\matrix(m)[matrix of math nodes,
row sep=4ex, column sep=4ex,
text height=1.5ex, text depth=0.25ex]
{&y&\\
x&&z\\};
\path[>=stealth,->,font=\scriptsize]
(m-1-2) edge[inner sep=0.75pt] node[above right]{$\psi$} (m-2-3)
(m-2-1) edge[inner sep=0.75pt] node[above left]{$\phi$} (m-1-2)
edge node[below]{$\chi$} (m-2-3);
\end{tikzpicture}
\end{equation*}
a $2$-simplex of $X\times_SK_{\alpha}^{\lhd}$ in which $\phi$ is $\pi_{f_{\alpha}}$-cocartesian and $\pi_{f_{\alpha}}(\psi)$ is an equivalence. Then the image of $\psi$ in $X$ is marked if and only if the image of $\chi$ is.
\end{itemize}
Finally, suppose
\[\rho\colon\fromto{(X,E)}{(S',M')}\]
a map of marked simplicial sets satisfying the following conditions.
\begin{itemize}
\item Any $2$-simplex of $X$ that lies over a scaled $2$-simplex of $S$ also lies over a scaled $2$-simplex of $S'$.
\item For any element $f_{\alpha}\colon\fromto{K_{\alpha}^{\lhd}}{S}$ of $P$ and any cocartesian section $s$ of $\pi_{f_{\alpha}}$, the composite
\[K_{\alpha}^{\lhd}\ \tikz[baseline]\draw[>=stealth,->,font=\scriptsize,inner sep=0.75pt](0,0.5ex)--node[above]{$s$}(0.5,0.5ex);\ X\times_SK_{\alpha}^{\lhd}\ \tikz[baseline]\draw[>=stealth,->](0,0.5ex)--(0.5,0.5ex);\ X\ \tikz[baseline]\draw[>=stealth,->,font=\scriptsize,inner sep=0.75pt](0,0.5ex)--node[above]{$\rho$}(0.5,0.5ex);\ S'\]
lies in $P'$.
\end{itemize}
Then the adjunction
\[\adjunct{\rho_!\circ\pi^{\star}}{s\Set^+_{/(S,M,T,P)}}{s\Set^+_{/(S',M',T',P')}}{\pi_{\star}\circ\rho^{\star}}\]
is a Quillen adjunction.
\end{thm}

Straight away, let us put this result to work.

\begin{cor}\label{cor:flatcocartRQF} Suppose $\pi\colon\fromto{X}{S}$ and $\rho\colon\fromto{X}{T}$ two functors of $\infty$-categories. Suppose $E\subset X_1$ a collection of marked edges on $X$, and assume the following.
\begin{itemize}
\item Any equivalence of $X$ is marked, and the marked edges are closed under composition.
\item The functor $\pi\colon\fromto{X}{S}$ is a flat locally cartesian fibration.
\item Suppose
\begin{equation*}
\begin{tikzpicture}[baseline]
\matrix(m)[matrix of math nodes,
row sep=4ex, column sep=4ex,
text height=1.5ex, text depth=0.25ex]
{&y&\\
x&&z\\};
\path[>=stealth,->,font=\scriptsize]
(m-1-2) edge[inner sep=0.75pt] node[above right]{$\psi$} (m-2-3)
(m-2-1) edge[inner sep=0.75pt] node[above left]{$\phi$} (m-1-2)
edge node[below]{$\chi$} (m-2-3);
\end{tikzpicture}
\end{equation*}
a $2$-simplex of $X$ in which $\pi(\phi)$ is an equivalence, and $\psi$ is locally $\pi$-cartesian. Then the edge $\phi$ is marked just in case $\chi$ is.
\end{itemize}
Then the functor $\pi_{\star}\circ\rho^{\star}$ is a right Quillen functor
\[\fromto{s\Set^+_{/T}}{s\Set^+_{/S}}\]
for the cocartesian model structure. In particular, if $\fromto{Y}{T}$ is a cocartesian fibration, then the map $r\colon\fromto{Z}{S}$ given by the universal property
\[\Mor_S(K,Z)\cong\Mor_T(K^{\flat}\times_{S^{\sharp}}(X,E),\leftnat{Y})\]
is a cocartesian fibration, and an edge of $Z$ is $r$-cocartesian just in case the map $\fromto{\Delta^1\times_SX}{Y}$ carries any edge whose projection to $X$ is marked to a $q$-cocartesian edge.
\end{cor}

We will speak of applying Cor. \ref{cor:flatcocartRQF} to the span
\[
S\ \tikz[baseline]\draw[>=stealth,<-,font=\scriptsize,inner sep=1pt](0,0.5ex)--node[above]{$\pi$}(0.5,0.5ex);\ (X,E)\ \tikz[baseline]\draw[>=stealth,->,font=\scriptsize,inner sep=1pt](0,0.5ex)--node[above]{$\rho$}(0.5,0.5ex);\ T
\]
to obtain a Quillen adjunction
\[
\adjunct{\rho_{!}\circ\pi^{\star}}{s\Set^+_{/S}}{s\Set^+_{/T}}{\pi_{\star}\circ\rho^{\star}}.
\]

\begin{exm}\label{exm:cartesianworkhorse} Suppose $p\colon\fromto{X}{S}$ is a cartesian fibration, and suppose $q\colon\fromto{Y}{S}$ is a cocartesian fibration. The construction of \ref{cor32213} gives a simplicial set $\widetilde{\underline{\Fun}}_S(X,Y)$ over $S$ given by the universal property
\begin{equation*}
\Mor_{S}(K,\widetilde{\underline{\Fun}}_S(X,Y))\cong\Mor_{S}(K\times_SX,Y).
\end{equation*}
By the previous result, the functor $\fromto{\widetilde{\underline{\Fun}}_S(X,Y)}{S}$ is a cocartesian fibration. Consequently, one finds that \cite[Cor. 3.2.2.13]{HTT} is a very elementary, very special case of the previous result.
\end{exm}

Another compelling example can be obtained as follows. Suppose $\phi\colon\fromto{A}{B}$ is a functor of $\infty$-categories. Of course if one has a cocartesian fibration $\fromto{Y}{B}$ classified by a functor $\YY$, then one may pull back to obtain the cocartesian fibration $\fromto{A\times_BY}{A}$ that is classified by the functor $\YY\circ\phi$. In the other direction, however, if one has a cocartesian fibration $\fromto{X}{A}$ that is classified by a functor $\XX$, then how might one write explicitly the cocartesian fibration $\fromto{Z}{B}$ that is classified by the right Kan extension? It turns out that this is just the sort of situation that Cor. \ref{cor:flatcocartRQF} can handle gracefully. We'll need a spot of notation:

\begin{ntn}\label{ntn:laxpb} If $f\colon\fromto{M}{S}$ and $g\colon\fromto{N}{S}$ are two maps of simplicial sets, then let us write
\[M\downarrow_SN\coloneq M\underset{\Fun(\Delta^{\{0\}},S)}{\times}\Fun(\Delta^1,S)\underset{\Fun(\Delta^{\{1\}},S)}{\times}N.\]
A vertex of $M\downarrow_SN$ is thus a vertex $x\in M_0$, a vertex $y\in N_0$, and an edge $\fromto{f(x)}{g(y)}$ in $S_1$. When $M$ and $N$ are $\infty$-categories, the simplicial set $M\downarrow_SN$ is a model for the lax pullback of $f$ along $g$.
\end{ntn}

\begin{prp}\label{exm:exceptionalQuillenAdjunction} If $\phi\colon A \to B$ is a functor of $\infty$-categories, and if $p\colon\fromto{X}{A}$ is a cocartesian fibration, then form the simplicial set $Y$ over $B$ given by the following universal property: for any map $\eta\colon\fromto{K}{Z}$, we demand a bijection
\[
\Mor_{/B}(K,Y)\cong\Mor_{/A}(K\downarrow_{B}A,X),
\]
functorial in $\eta$. Then the projection $q\colon\fromto{Y}{B}$ is a cocartesian fibration, and if $p$ is classified by a functor $\XX$, then $q$ is classified by the right Kan extension of $\XX$ along $\phi$.
\begin{proof} As we have seen, we need conditions on $\phi$ to ensure that the adjunction
\[ \adjunct{\phi^\ast}{s\Set^+_{/B}}{s\Set^+_{/A}}{\phi_\ast} \]
is Quillen. However, we can rectify this situation by applying Cor. \ref{cor:flatcocartRQF} to the span
\[
B\ \tikz[baseline]\draw[>=stealth,<-,font=\scriptsize,inner sep=1pt](0,0.5ex)--node[above]{$\ev_0 \circ \pr_{\mathscr{O}(B)}$}(1,0.5ex);\ (\mathscr{O}(B) \times_{\ev_1,B,\phi} A)^\sharp\ \tikz[baseline]\draw[>=stealth,->,font=\scriptsize,inner sep=1pt](0,0.5ex)--node[above]{$\pr_A$}(0.5,0.5ex);\ A
\]
in view of the following two facts, which we leave to the reader to verify:
\begin{itemize}
\item the functor $\ev_0 \circ \text{pr}_{\mathscr{O}(B)}$ is a cartesian fibration;
\item for all $C \to B$ fibrant in $s\Set^+_{/B}$, the identity section $\fromto{B}{\mathscr{O}(B)}$ induces a functor
\[\fromto{C \times_B A^\sharp}{C \times_B \mathscr{O}(B)^\sharp \times_B A^\sharp}\]
which is a cocartesian equivalence in $s\Set^+_{/A}$.
\end{itemize}

We thereby obtain a Quillen adjunction
\[ \adjunct{(\text{pr}_A)_! \circ (\ev_0 \circ \text{pr}_{\mathscr{O}(B)})^\ast}{s\Set^+_{/B}}{s\Set^+_{/A}}{(\ev_0 \circ \text{pr}_{\mathscr{O}(B)})_\ast \circ (\text{pr}_A)^\ast} \]
which models the adjunction of $\infty$-categories
\[ \adjunct{\phi^\ast}{\Fun(B,\Cat_\infty)}{\Fun(A,\Cat_\infty)}{\phi_\ast}, \]
as desired.
\end{proof}
\end{prp}

This result illustrates the basic principle that replacing strict pullbacks by lax pullbacks allows one to make homotopically meaningful constructions. Of course, any suitable theory of $\infty$-categories would allow one to produce the right adjoint $\phi_\ast$ less explicitly, by recourse to the adjoint functor theorem. This example rather illustrates, yet again, how the model of quasicategories allows one to control common categorical constructions in an utterly transparent manner.

Let us end our study with some simple but gratifying observations concerning the interaction of Th. \ref{thm:functorialityofcatpatterns} with compositions and homotopy equivalences of spans. First the compositions:
\begin{lem} \label{spanCompositionLem} Suppose we have spans of marked simplicial sets
\[
C_0\ \tikz[baseline]\draw[>=stealth,<-,font=\scriptsize,inner sep=1pt](0,0.5ex)--node[above]{$\pi_0$}(0.5,0.5ex);\ D_0\ \tikz[baseline]\draw[>=stealth,->,font=\scriptsize,inner sep=1pt](0,0.5ex)--node[above]{$\rho_0$}(0.5,0.5ex);\ C_1
\]
and
\[
C_1\ \tikz[baseline]\draw[>=stealth,<-,font=\scriptsize,inner sep=1pt](0,0.5ex)--node[above]{$\pi_1$}(0.5,0.5ex);\ D_1\ \tikz[baseline]\draw[>=stealth,->,font=\scriptsize,inner sep=1pt](0,0.5ex)--node[above]{$\rho_1$}(0.5,0.5ex);\ C_2
\]
and categorical patterns $(M_i,T_i,\emptyset)$ on $C_i$ such that the hypotheses of Th. \ref{thm:functorialityofcatpatterns} are satisfied for each span. Then the hypotheses of Th. \ref{thm:functorialityofcatpatterns} are satisfied for the span 
\[
D_0\ \tikz[baseline]\draw[>=stealth,<-,font=\scriptsize,inner sep=1pt](0,0.5ex)--node[above]{$\pr_{0}$}(0.5,0.5ex);\ D_0\times_{C_1}D_1\ \tikz[baseline]\draw[>=stealth,->,font=\scriptsize,inner sep=1pt](0,0.5ex)--node[above]{$\pr_{1}$}(0.5,0.5ex);\ D_1
\]
with the categorical patterns $(M_{D_i}, \pi_i^{-1}(T_i),\emptyset)$ on each $D_i$. Consequently, we obtain a Quillen adjunction
\[
\adjunct{(\rho_1\circ\pr_1)_!\circ(\pi_0\circ\pr_0)^{\star}}{s\Set^+_{/(C_0,M_0,T_0,\varnothing)}}{s\Set^+_{/(C_2,M_2,T_2,\varnothing)}}{(\pi_0\circ\pr_0)_{\star}\circ(\rho_1\circ\pr_1)^{\star}},
\]
which is the composite of the Quillen adjunction from $s\Set^+_{/(C_0,M_0,T_0,\varnothing)}$ to $s\Set^+_{/(C_1,M_1,T_1,\varnothing)}$ with the one from $s\Set^+_{/(C_1,M_1,T_1,\varnothing)}$ to $s\Set^+_{/(C_2,M_2,T_2,\varnothing)}$.
\begin{proof} The proof is by inspection. However, one should beware that the ``long'' span
\[
C_0\ot\ D_0\times_{C_1}D_1\to D_1
\]
can fail to satisfy the hypotheses of Th. \ref{thm:functorialityofcatpatterns}, because the composition of locally cartesian fibrations may fail to be again be locally cartesian; this explains the roundabout formulation of the statement. Finally, observe that if we employ the base-change isomorphism $\rho_0^\ast \pi_{1,\ast} \cong \pr_{0,\ast} \circ \pr_1^\ast$, then we obtain our Quillen adjunction as the composite of the two given Quillen adjunctions.
\end{proof}
\end{lem}

\noindent Now let us see that Th. \ref{thm:functorialityofcatpatterns} is compatible with homotopy equivalences of spans:
\begin{lem} \label{spanFunctorialityLem} Suppose a morphism of spans of marked simplicial sets
\[ \begin{tikzpicture}[baseline]
\matrix(m)[matrix of math nodes,
row sep=6ex, column sep=4ex,
text height=1.5ex, text depth=0.25ex]
 {  & K &  \\
  C & L & C' \\ };
\path[>=stealth,->,font=\scriptsize]
(m-1-2) edge[inner sep=0.75pt] node[left]{$f$} (m-2-2)
edge node[above left]{$\pi$} (m-2-1)
edge node[above right]{$\rho$} (m-2-3)
(m-2-2) edge node[below]{$\pi'$} (m-2-1)
edge node[below]{$\rho'$} (m-2-3);
\end{tikzpicture} \]
where $\rho_! \pi^\ast$ and $(\rho')_! (\pi')^\ast$ are left Quillen with respect to the model structures given by categorical patterns $\mathfrak{P}_C$ and $\mathfrak{P}_{C'}$ on $C$ and $C'$. Suppose moreover that $f$ is a homotopy equivalence in $s\Set^+_{/\mathfrak{P}_{C'}}$ -- i.e., suppose that there exists a homotopy inverse $g$ and homotopies
\[h\colon \id \simeq g \circ f\text{\quad and\quad}k\colon \id \simeq f \circ g.\]
Then the natural transformation $\rho_! \pi^\ast \to (\rho')_! (\pi')^\ast$ induced by $f$ is a weak equivalence on all objects, and, consequently, the adjoint natural transformation $(\pi')_\ast (\rho')^\ast \to \pi_\ast \rho^\ast$ is a weak equivalence on all fibrant objects.
\begin{proof} The homotopies $h$ and $k$ pull back to show that for all $X \to C$, the map
\[\id_X \times_C f\colon X \times_C K \to X \times_C L\]
is a homotopy equivalence with inverse $\id_X \times_C g$. The last statement now follows from \cite[1.4.4(b)]{Hovey}.
\end{proof}
\end{lem}

%-------------------------------------------------------------------%
%-------------------------------------------------------------------%
%-------------------------------------------------------------------%

\bibliographystyle{amsplain}
\bibliography{Gcats}

\end{document}